\input amstex
\documentstyle{amsppt}
\tolerance=2000 
\define\m1{^{-1}}
\define\ov1{\overline}
\def\gp#1{\langle#1\rangle}
\def\ul2#1{\underline{\underline{#1}}}
\def\lemma#1{ Lemma #1}
\def\proof{{\it Proof.} \hskip 4pt}

\TagsOnRight

\topmatter
\title 
Group algebras whose involutory units commute
\endtitle

\dedicatory
Dedicated to the memory of Professor I.\ I.\ Khripta
\enddedicatory

\author
Victor Bovdi$^1$ and Michael Dokuchaev$^2$
\endauthor
\affil
${}^1$Institute of Mathematics and Informatics, Debrecen University, Hungary 
\\
$vbovdi\@math.klte.hu$
\\
${}^2$Departamento de Matem\'atica, Universidade de S\~ao Paulo, Brasil
\\
$dokucha\@ime.usp.br$
\endaffil

\leftheadtext\nofrills{Victor Bovdi and Michael Dokuchaev}
\rightheadtext\nofrills{ Group algebras whose involutory units commute
}

\abstract
Let $K$ be a field  of characteristic $2$ and $G$ a nonabelian 
locally finite $2$-group.
Let $V(KG)$be  the group of units with augmentation $1$ 
in the group algebra
$KG$. An explicit list of groups is given, and it is proved 
that all
involutions in $V(KG)$ commute with each other if and only 
if $G$ is isomorphic
to one of the groups on this list. In particular, this 
property depends
only on $G$ and not at all on $K$.
\endabstract

\thanks 
The author's research was supported  by OTKA  No.T 025029; No. T 029132  
and by FAPESP Brazil (proc.\ 97/05920-6).
The second author was partially supported by CAPES Brazil 
(proc. AEX0871/96-6) and partially by CNPq Brazil (proc.\
301115195-8)
\endthanks 

\endtopmatter
\document

\subhead
Introduction
\endsubhead

Let  $KG$ be the group algebra of a locally finite $p$-group $G$ 
over a field  $K$ of characteristic $p$. Then the normalized unit 
group
$$
 V(KG)=\Bigl\{\,\sum_{g\in G} \alpha_g g\in KG\Bigm|
 \sum_{g\in G} \alpha_g=1\,\Bigr\}
$$
is a locally finite $p$-group.

An interesting way to study  $V(KG)$ is to construct embeddings 
of important groups into it. 
D.B. Coleman and D.S. Passman  \cite{4} have proved that if $G$ is 
nonabelian then the wreath product of two cyclic groups of order $p$ is
involved in $V(KG)$. Larger wreath products have been constructed by  A. Mann and  A. Shalev in \cite{7}, \cite{8}, \cite{9} and  \cite{10}. In this paper we answer the question of when dihedral groups 
cannot be embedded into  $V(KG)$. Clearly, only the case $p=2$ 
has to be considered. At the same time we obtain the list of the 
locally finite $2$-groups $G$ such that $V(KG)$  does not contain 
a subgroup isomorphic to a wreath product of two groups (for the 
case of odd $p$ see \cite{1}). 

We remark, that in the present work the method of proof developed
 in \cite{3} plays an essential role.   

Let $C_{2^n}$, $C_{2^\infty}$ and  
$Q_8$  be the  cyclic group of order $2^n$, the   quasicyclic group 
of type $2^\infty$  and the  quaternion group of order $8$,  
respectively.  An {\it involution } is a group element of order $2$. 
For any $a,b\in G$  we denote 
$[a,b]=a\m1b\m1ab$, $a^b=b\m1ab$.

Our main result is the following:

\proclaim {Theorem }
{\it
Let  $K$ be a field of characteristic $2$, 
and $G$  a locally  finite nonabelian $2$-group. Then  
all involutions of $V(KG)$ commute 
if and only if $G$ is one of the following groups:
\itemitem{\rm(i)} 
$S_{n,m}=\gp{a,b\mid a^{2^n}=b^{2^m}=1, 
a^b=a^{1+2^{n-1}}}$ with $n,m\geq 2$,
or  $Q_8$;
\itemitem{\rm(ii)} $Q_8\times C_{2^n}$ or 
$Q_8\times C_{2^\infty}$; 
\itemitem{\rm(iii)} the semidirect product of 
the cyclic group  $\gp{d \mid d^{2^n}=1}$  with the 
quaternion group
$\gp{\,a,b\mid a^4=1, a^2=b^2=[a,b]}$ such that
$[a,d]=d^{2^{n-1}}$ and $[b,d]=1$;
\itemitem{\rm(iv)} $H_{32}=\bigl\langle\,x,y,u \bigm|\
x^4=y^4=1,\  x^2=[y,x],
y^2=u^2=[u,x], \ x^2y^2=[u,y]\,\bigr\rangle$.

}
\endproclaim

For a 
$2$-group $G$ we denote by $\Omega(G)$ the subgroup 
generated by all elements of order $2$ of $G$.
We write  
$\Omega(V)=\Omega\big(V(KG)\big)$. Clearly, $\Omega(G)$ is a normal 
subgroup of $G$. As usually,
$\exp(G)$,  $C_G\big(\gp{a,b}\big)$
denote the exponent of a group $G$ and the centralizer of the 
subgroup $\gp{a,b}$ in $G$, respectively.  
Let
$$
{\Cal L}(KG)=\gp{xy-yx\mid x,y\in KG}
$$ 
be the  commutator  subspace of $KG$.

For an element $g$  of a finite order $|g|$ in a group $G$, 
let $\ov1{g}$ denote  the sum (in $KG$) of the distinct 
powers of $g$:
$$
\ov1{g}=\sum_{i=0}^{|g|-1}g^i.
$$  
For an  arbitrary element $x=\sum_{g\in G}\alpha_gg\in KG$ 
we put  $\chi(x)=\sum_{g\in G}\alpha_g\in K$.

It will be convenient to have a short
temporary name for the locally finite $2$-groups $G$ such 
that all elements of order $2$ of $V(KG)$ form an abelian  
subgroup.  Let us call the 
groups $G$ with this property 
{\it good.}

\subhead
1. Preliminary results
\endsubhead

\proclaim
{\lemma 1}
{\it
Let $G$ be a finite nonabelian good  group. Then all 
involutions of $G$ are central, 
$G^\prime\subseteq \Omega(G)$ and
either $G$ is  the quaternion group of order $8$
or $\Omega(G)$ is a direct product of two  cyclic groups.
}
\endproclaim

\proof
Let $G$ be a nonabelian good  group. Clearly
$\Omega(V)\cap G=\Omega(G)$ is a normal abelian
subgroup of $G$.

Suppose that $|\Omega(G)|=2$. Then by Theorem 3.8.2 in \cite{5} 
we have that 
$$
G=\gp{a,b\mid a^{2^m}=1, b^2=a^{2^{m-1}}, a^b=a\m1}
$$
with $ m\geq 2$ 
is a generalized quaternion $2$-group.  If $|G|>8$, then we choose
$c\in\gp{a}$ of order $8$.  Then $1+(1+c^2)(c+b)$ and
$1+(1+c^2)(c+cb)$ are noncommuting involutions of $V(KG)$, which
is impossible. Therefore, $G$ is the quaternion group of order $8$.

Let $|\Omega(G)|>2$. The normal subgroup
$\Omega(G)$ contains a central element $a$ in $G$ and $x=1+(a+1)g$
is an involution for any $g\in G$. If $b\in \Omega(G)$, then $bx=xb$ and
$(a+1)(1+[b,g])=0$. 
It  implies that $[g,b]=1$ or $[g,b]=a$. Let $[g,b]=a$ and $|g|=2^t$. 
Then $z=1+(g+1)\ov1{g^2}$ is an involution and $zb=bz$. From this 
we get that  $a \in \gp{g^2}$ and since $a\in \Omega(G)$ we obtain  that 
$a=g^{2^{t-1}}$. Clearly, if 
$t=2$, then $\gp{b,g}$ is the dihedral group of order $8$, which 
is impossible. If $t>2$, then $x=1+g(1+g^{2^{t-2}})(1+b)$ is 
an involution which does  not commute with $b$, which is a contradiction. 
Therefore, $[g,b]=1$ and  $\Omega(G)$ is a central subgroup
of $G$.  Let $a$ and $b$ be arbitary elements of 
$\Omega(G)$, $g,h\in G$ and $[g,h]\not=1$. Then $x=1+(a+1)g$ 
and $y=1+(b+1)h$ are involutions and  $[x,y]=1$. 
From this we conclude that the commutator subgroup  $G^\prime$ is
a subgroup of $\gp{a, b}$ in $\Omega(G)$ and $\exp(G^\prime)=2$.

Now, let $|\Omega(G)|\geq 8$ and let $a$, $b$, $c$ be linearly 
independent elements of $\Omega(G)$.  Then 
by the above reasoning 
we
 have $[g,h]\in
\gp{a,b}\cap\gp{a,c}\cap\gp{b,c}=1$, which is impossible. Therefore, 
$|\Omega(G)|=4$ and  $G'\subseteq
\Omega(G)$.

\rightline{\text{\qed}}

\proclaim
{\lemma 2}
{\it
A two-generator finite nonabelian  group  is good if and only if it
is either the quaternion group of order $8$ or
$$
S_{n,m}=\gp{a,b\mid a^{2^n}=b^{2^m}=1,
a^b=a^{1+2^{n-1}}}\tag1
$$
with $n,m\geq2$.
}
\endproclaim
\proof
Suppose that $G$ is not the quaternion group of order $8$. By
Lemma 1 we have that  $\Omega(G)$ is a direct product of two  cyclic groups 
and $G^\prime\subseteq \Omega(G)$. Then the Frattini subgroup
$\Phi(G)=\{ g^2 \mid g\in G\}$ is  central and by Theorem
3.3.15 in \cite{5}, $|G/\Phi(G)|=4$. Since $\Phi(G)$ is a
subgroup of the centre $\zeta(G)$ and the factor group 
$G/\zeta(G)$ can not be cyclic, this implies $\Phi(G)=\zeta(G)$.

It is easy to verify that a two-generator good group $G$ is
metacyclic. Indeed, every maximal subgroup $M$ of $G$ is abelian
and normal in $G$, because $\Phi(G)=\zeta(G)\subset M$ and
$|M/\zeta(G)|=2$. Clearly, $\Omega(M)\subseteq \Omega(G)$ and
in case $|\Omega(M)|=2$ the subgroup $M$ is cyclic and
we conclude that $G$ is  metacyclic.

Now let $|\Omega(M)|=4$ for every maximal subgroup $M$ of $G$. 
Then $G$ and $M$ are two-generator groups. It is easy to see that  all 
such groups of order $16$ are metacyclic. If $|G|=2^n$ 
($n\geq 5$) then   $G$ and all maximal subgroups of $G$ are 
two generator groups and by Theorem 3.11.13 in \cite{5} $G$ is 
metacyclic too. Since $G^\prime\subseteq \Omega(G)$
by Lemma 1, it follows from  Theorem 3.11.2 in \cite{5},  $G$ is
defined by
$$
 G=\gp{a,b\mid a^{2^n}=1, b^{2^m}=a^{2^l}, 
 a^b=a^{1+2^{n-1}}}
$$
with  $n,m\geq 2$.

Suppose that the generators $a$ and  $b$ are choosen with minimal   order of $b$. We want to show  that $b^{2^m} = 1$. Suppose by contradiction that the order of $b$ is bigger than $2^m$. Then $m \leq l$. Indeed, if $m > l$  then $b^{2^{m-l}} \in \zeta(G)$ and        $ (ab^{-2^{m-l}})^{2^l} = 1$, so we may take generators $a_1 = b $ and $ b_1 = ab^{-2^{m-l}}$, but the order of $b_1$ is less than  $ 2^m$, a contradicition with our choice of generators.
Thus, $m \leq l$.  If $m < l$, then $a^{2^{l-m}}$ is central and, taking $a_1 = a, b_1 = a^{-2^{l-m}} b$, we have that 
$(a^{-2^{l-m}} b )^{2^m} = 1$ and the order of $b_1$ is less then the order of $b$, which is not allowed. Now,  the case  $l=m > 1$ is also  impossible, since we can take generators $a_1 = a, b_1 = a^{-1}b$ with $b_1^{2^m}=1$. In the last case, when $l=m=1$, we have $n \geq 3$ as $G$ is not quaternion of order $8$. Hence, $a^{2^{n-2}} \in \zeta(G)$ and $a^{2^{n-2}-2}ab$ is a non-central element of order $2$, which gives the final contradicition.

\rightline{\text{\qed}}

\proclaim
{\lemma 3} {\rm (E.~A.~O'Brien, see Lemma 4.1 in \cite{3})}
{\it
The groups $H$ of order dividing $128$ in which $\Phi(H)$
and $\Omega(H)$ are equal, central, and of order $4$, are
precisely the following: $C_4\times C_4$, $C_4\rtimes C_4$,
$C_4\rtimes Q_8$, $Q_8\times C_4$, $Q_8\times Q_8$,
the central product of the group
$S_{2,2}=\gp{\,a,b\mid a^4=b^4=1,\ a^2=[b,a]\,}$ with a
quaternion group of order $8$, the nontrivial element
common to the two central factors being $a^2b^2$,
$$
\eqalign{
 H_{245}=\bigl\langle\,x,y,u,v\bigm|\ &x^4=y^4=[v,u]=1,\cr
&x^2=v^2=[y,x]=[v,y],\cr  &y^2=u^2=[u,x],\cr
&x^2y^2=[u,y]=[v,x]\,\bigr\rangle\cr}
$$
and the  groups named in parts {\rm(iii), (iv)} of Theorem.

}
\endproclaim

\rightline{\text{\qed}}

 The  group $H_{245}$, is one of the
two Suzuki $2$-groups  of order $64$.

\subhead
2. Proof of the `only if' part of the theorem
\endsubhead

Let $G$ be a finite nonabelian good  group, so, $\Omega(V)$ is abelian. 

By Lemma 1   all involutions of $G$ are central, $G^{\prime}\subseteq 
\Omega(G)$  and
$\Omega(G)$ is either a group of order $2$ or a direct 
product of two  cyclic groups.  Clearly, 
$\Phi(G)\subseteq \zeta(G)$ and if $|\Omega(G)|=2$ then 
by Lemma 1 $G$ is a quaternion group of order $8$. Thus, 
we can suppose that $|\Omega(G)|=4$.

First let $\Phi(G)$ be  cyclic. Since all 
involutions are central, by Theorem 2 in \cite{2} $G$ is 
the direct product of a group of order $2$ and the generalized 
quaternion group of order $2^{n+1}$. By Lemma 1 $G$ is a 
Hamiltonian $2$-group of order $16$.

We may suppose that $\Phi(G)$ is the direct product of two cyclic groups.
Let the exponent of $G$ be  $4$. Then 
$\Phi(G)=\Omega(G)$ and  by a result of N.~Blackburn 
(Theorem VIII.5.4 in \cite{6}), $|G|\leq |\Omega(G)|^3$. 
Therefore, the order of $G$ divides $64$.  Then   by 
O'Brien's Lemma,  $G$ is precisely one 
of the following types: $C_4\rtimes C_4$, $C_4\rtimes Q_8$,
$Q_8\times C_4$, the groups named in parts {\rm(iii), (iv)} 
of Theorem and  $Q_8\times Q_8$, $H_{245}$, the central product of $Q_8$ 
and  $S_{2,2}$ with common  $a^2b^2$.

Now we shall find noncommuting involutions $z_1,z_2$ in 
$V(KG)$ if $G$ is one of the last three groups listed above.

Let $G$ be  the  central product of the group
$$
S_{2,2}=\gp{\,a,b\mid a^4=b^4=1,\ a^2=[b,a]\,}
$$ 
with the quaternion group of order $8$, the nontrivial element 
common to the two central factors being $a^2b^2$. Then 
$$
\eqalign{
G\cong \gp{a,b,d,f\mid &a^4=d^4=1, b^2=a^2=[a,b], f^2=d^2=[d,f],\cr
&[a,d]=[b,d]=[b,f]=1, [a,f]=a^2}}
$$
and we put 
$z_1=1+d^2a+b+a^3d+bd+f+abf+df+abdf$ and $z_2=1+b(1+d^2)$.
If $G\cong  H_{245}$ then  
$$
\eqalign{
 H_{245}\cong\bigl\langle\,
a,b,d,f\bigm|\ &a^4=b^4=1, b^2=d^2=a^2, [a,b]=1,\cr  
&[a,d]=[b,f]=[d,f]=b^2, [b,d]=a^2, [a,f]=a^2b^2\,\bigr\rangle\cr}
$$
and put 
$z_1=1+a+ab+d+a^2bd+f+bf+ab^2df+a^3b^3df$ and $z_2=1+(b+b\m1)$.
Now, let $G$ be a  direct product of two quaternion groups 
$\gp{a,b}$ and $\gp{c,d}$ of order $8$. Then we put 
$z_1= 1+a+bc^2+c+abc+a^2d+abd+acd+bcd$   and $z_2=1+b(1+c^2)$. 

It is easy to verify that in all three cases $z_1^2=z_2^2=1$, 
$z_1z_2\not=z_2z_1$.

Now, let the exponent of  $G$ be greater than $4$.
Using Lemma 2, we conclude that $G$ contains a  two-generator
nonabelian subgroup $H$ which is either $Q_8$ or $S_{n,m}$.

We wish to prove that if $\exp(G)> 4$ and $G=H\cdot C_G(H)$  for every 
two-generator nonabelian subgroup $H$, then 
$$
G=Q_8\times\gp{d\mid d^{2^n}=1, n>1}.
$$
First, let $H=Q_8=\gp{ a,b}$ be a  quaternion subgroup of order 
$8$ of $G$. Then $G=Q_8\cdot C_G(Q_8)$ and  $C_G(Q_8)$ does not
contain an element $c$ of order $4$ with the property $c^2=a^2$,
because $ac$ would be  a noncentral involution of $G$, which is impossible.
If $C_G(Q_8)$ is abelian and $|\Omega\big(C_G(Q_8)\big)|=4$ then $C_G(Q_8)$
is the direct product of $\gp{a^2}$ and
$\gp{d\mid d^{2^n}=1}$, with  ${n>1}$, and $G=Q_8\times \gp{d}$.

We can suppose that $C_G(Q_8)$ is nonabelian and does not contain an
element $u$ such that $u^2=a^2$. Since $\exp\big(C_G(Q_8)\big)>4$, there
always exists a subgroup 
$$
S_{n,m}=\gp{c,d\mid c^{2^n}=d^{2^m}=1, 
c^d=c^{1+2^{n-1}}}
$$ 
of $C_G(Q_8)$ which is of exponent greater than $4$.
Then $\zeta(S_{n,m})=\gp{c^2,d^2}$ and since   
$\exp(S_{n,m})>4$, one of the generators $c$ or $d$ has order 
greater than $4$. Therefore, any $u\in \Omega\big(\zeta(S_{n,m})\big)$ is 
the  square of one of the elements from $S_{n,m}$. 
Thus,  since we assume that $C_G(Q_8)$ does not contain an element 
of order $4$ whose square is $a^2$, we get  that $S_{n,m}\cap Q_8=1$ and
$S_{n,m}\times Q_8$ is a subgroup of $G$ with the property
$|\Omega(S_{n,m}\times Q_8)|=8$, which is impossible.

Now let $H=S_{n,m}=\gp{a,b\mid a^{2^n}=b^{2^m}=1,
a^b=a^{1+2^{n-1}}}$ with ${ n,m\geq 2}$
be a subgroup of $G$.  Then $G=S_{n,m}\cdot C_G(S_{n,m})$.
Since  $|\Omega(G)|=4$,  we can 
choose $d\in C_G(S_{n,m})$ such that $d\not\in S_{n,m}$ but
$d^2\in S_{n,m}$. Then $d^2= a^{2i}b^{2j}$. If $i$ or  $j$
is even then  $d\m1a^ib^j\in \Omega(G)=\Omega(S_{n,m})$ and
$d\in S_{n,m}$, which is impossible.
If $i$ and $j$ are odd and $n>2$ then $d\m1a^{i+2^{n-2}}b^j\in
\Omega(G)=\Omega(S_{n,m})$ which  is  a contradiction.
Therefore, $n=2$ and $\gp{a^i,d\m1b^j}$ is a
quaternion subgroup and by assumption $G=Q_8\cdot C_G(Q_8)$. We
obtained the previous case.

It is easy to check that  if the commutator subgroup $G^{\prime}$ 
is of order $2$, then   $G=H\cdot C_G(H)$  for  every  
two-generator nonabelian subgroup of $G$. Indeed, 
the equalities $G^{\prime}=H^{\prime}=\gp{c}$  imply 
that $H=\gp{a,b}$ is a normal 
subgroup of $G$. Let 
$[a,b]=c$, $[a,g]=c^k$, $[b,g]=c^l$, where 
$0\leq k,l\leq 1$,  $g\in G$. Then at least one of the elements $g$, 
$ag$, $bg$, $abg$ belongs to $C_G(H)$ and $g\in H\cdot C_G(H)$.
Therefore, $G=H\cdot C_G(H)$.

It follows that we can suppose that  $\exp(G)>4$, the commutator 
subgroup $G'=\Omega(G)$
has order $4$ and $G$ contains a two-generator nonabelian subgroup
$L$ such that $G\not=L\cdot C_G(L)$.

Let $L=\gp{b,d\mid [b,d]\not=1}$ and $a\in G\setminus (L\cdot C_G(L))$. Then
$[a,d]$ or $[b,a]$ is not equal to $[b,d]$. Indeed, in the contrary
case, from $[a,b]=[d,b]$ and $[a,d]=[b,d]$ we get  
$bda\in C_G(L)$ and $a\in L\cdot C_G(L)$ which is
impossible.

Now we want to prove that we  can choose 
$a\in G\setminus (L\cdot C_G(L))$ and 
$b, d\in L$ such that $\gp{b,d}=L$,  $[a,b]=1$
with the following property:
$$
[a,d]\not=[b,d], [a,d]\not=1, [b,d]\not=1. \tag2
$$ 

If $[a,d]=1$ then we can put $a'=a, b'=d$ and $d'=b$.

We consider  the following cases:

Case 1. Let $[b,d]=[b,a]\not=[a,d]$. Then $[b,ad]=1$ and we put $a'=ad$, 
$b'=b$ and $d'=d$. If $[a',d']=1$ then $ad\in C_G(L)$ which implies that 
$a\in L\cdot C_G(L)$, a contradiction. 

Case 2. Let $[b,d]=[a,d]\not=[a,b]$. Then $[ab,d]=1$ and  put $a'=ab$, 
$b'=d$ and $d'=b$. If $[a',d']=1$  then  $ab\in C_G(L)$ which gives 
a contradiction again.

Case 3. Let $[a,b]\not=[a,d]\not=[b,d]\not=[a,b]$. 
Suppose that all these commutators are not trivial. 
Since $|\Omega(G)|=4$,
one of these commutators  equals to the product of two others and 
$$
[ab,bd]=[a,b]\cdot [b,d]\cdot [a,d]=1.
$$
Put $a'=ab$, $b'=bd$ and $d'=d$. 

In all what follows we suppose that $L=\gp{b,d}$ and 
$a\in G\setminus (L\cdot C_G(L))$ such that $[a,b]=1$ and   (2) 
is satisfied.  

It is easy to see that if $\gp{a,b}=\gp{u}$ is cyclic then from 
$[a,d]\not=1\not=[b,d]$  we have $a=u^{2k+1}$ and $b=u^{2t+1}$ 
for some $k,t\in \Bbb N$, because the squares of  all elements 
in $G$ are  central. Then $ab\in C_G(L)\subseteq L\cdot C_G(L)$ 
and $a\in L\cdot C_G(L)$, which is impossible. Therefore, 
$\gp{a,b}$ is not cyclic.

Consider, $W=\gp{a,b,d}$. Then the commutator subgroup of $W$ has 
order $4$ and 
$$
H=\gp{a,b}=\gp{a_1}\times \gp{b_1}.
$$ 
Clearly $W=\gp{a_1,b_1,d}$ and $|W'|=4$. It is easy to see that 
$a_1$ and  $b_1$ can be choosen such that  condition   (2) is satisfied and  
$\gp{a_1}\cap \gp{b_1}=1$.

Let $a,b,d\in G$ with the property   (2), $[a,b]=1$ and 
$\gp{a}\cap\gp{b}=1$. Put 
$H=\gp{a \mid a^{2^n}=1}\times \gp{b \mid b^{2^m}=1}$ and 
$W=\gp{a,b,d}$.
Then 
$$
G'=\Omega(G)=\Omega(H)=\gp{a^{2^{n-1}}}\times \gp{b^{2^{m-1}}}
$$ 
and $H$ is a normal subgroup of $G$.

First, we shall  prove that $g^2\in H$ for  every $g\in G\setminus H$.
There exists $c=g^{2^{k-1}}$ such that $c\not\in H$ and $c^2\in H$.
If $k>1$, then $c\in \Phi(G)\subset \zeta(G)$ and we obtain that 
$c^2\in \zeta(G)\cap H$ and $c^2=a^{2i}b^{2j}$. Thus, $(c\m1a^ib^j)^2=1$ and 
$c\m1a^ib^j\in \Omega(G)=\Omega(H)$.
It  implies that $c\in H$, which is impossible. Therefore, $k=1$
and $g^2=a^{2i}b^{2j}$ for some $i$ and $j$, and we have shown that 
$g^2\in H$ for all $g\notin H$. 

We have 
$[g,a^ib^j]=a^{s2^{n-1}}b^{r2^{m-1}}$ for some $r,s\in \{0,1\}$.
It is easy to see that the case $n>2$ and $m>2$ is impossible. 
Indeed, if $n>2$ and $m>2$ then 
$a^{i+s2^{n-2}}\cdot b^{j+r2^{m-2}}\cdot g\m1\in \Omega(G)=\Omega(H)$ and 
$g\in H$, which is a contradiction.

Since $\exp(G)>4$ and for any $g\in G$, $g^2=a^{2i}b^{2j}$ 
for some $i,j$ it follows that $\exp(H)>4$. Thus 
 we may suppose that $n>2$ 
and $b$ is an element of order $4$.

Now we describe the group $W=\gp{a,b,d}$ and we distinguish
a number of cases according to the form of the element $d^2$.

Case 1. Let $d^2=a^{2i}$. Since $d\not\in H$ and $a^id\m1$ is 
not of  order $2$,
we have that $i=2k+1$ is odd. Then as $(a^{i+2^{n-2}}d\m1)^2\not=1$ 
we have 
$$
\big(a^{i+2^{n-2}}d\m1\big)^2=[a^i,d]\cdot a^{2^{n-1}}
=[a^{2k+1},d]\cdot a^{2^{n-1}}=
[a,d]\cdot a^{2^{n-1}}\not=1.
$$
It follows that $[a,d]\not= a^{2^{n-1}}$ and
$[a,d]=a^{s2^{n-1}}b^2$, $s\in\{0, 1\}$, and by property   (2) 
$[b,d]=a^{(1+s)2^{n-1}}b^2$ or  $[b,d]=a^{2^{n-1}}$. 
If $[b,d]=a^{2^{n-1}}$ then 
$a^{i+(1+s)2^{n-2}}bd\m1$ has order $2$ and does not belong to  
$\Omega(H)$ and this case is impossible.

Let $[a,d]=a^{s2^{n-1}}b^2$ and $[b,d]=a^{(1+s)2^{n-1}}b^2$.
Then
$$
W=\gp{a,b,d}=\gp{a^{i+2^{n-2}}d\m1, a^{(1+s)2^{n-2}}b, ab}
$$
and $\gp{a^{i+2^{n-2}}d\m1, a^{(1+s)2^{n-2}}b}$ is the  
quaternion subgroup of order $8$. Moreover, $ab$ has 
order $2^n$ and 
$[a^{i+2^{n-2}}d\m1,ab]=a^{2^{n-1}}=(ab)^{2^{n-1}}$. 
This shows  that $W$ satisfies (\rm{iii}) of the  Theorem.
Observe that $\gp{a^{(1+s)2^{n-2}}b, ab}=H$.

Case 2. Let $d^2=b^2$. As before,
$$
\big(a^{2^{n-2}}bd\m1\big)^2=a^{2^{n-1}}[b,d]\not=1.
$$
Therefore, $[b,d]\not= a^{2^{n-1}}$, and we obtain
$[b,d]=a^{s2^{n-1}}\cdot b^2$, and by property 
  (2) $[a,d]=a^{(1+sr)2^{n-1}}\cdot
b^{2r}$, where $r,s\in\{0,1\}$. It is easy to see that
$\gp{bd\m1, a^{s2^{n-2}}b}$ is the quaternion
group of order $8$ and $W=\gp{bd\m1, a^{s2^{n-2}}b, ab^r}$
is defined as in case $1$ and satisfies condition (iii) of the 
Theorem.
Moreover, $\gp{a^{s2^{n-2}}b,ab^r}=H$.

Case 3. Let $d^2=a^{2i}b^2$. If $i$ is even then $a^i\in \zeta(G)$
and $(da^{-i})^2=b^2$. 
Then $W=\gp{a,b, da^{-i}}$ and if we replace  $d$ by $d'=da^{-i}$ we obtain 
${d'}^2=b^2$ which is  case 2.

Let now $d^2=a^{2i}b^{2}$ and $i$ is odd. If 
$[a,d]=a^{s2^{n-1}}\cdot b^2$, then $a^{i+s2^{n-2}}\cdot d\m1\in
\Omega(H)$ and $d\in H$, which is impossible. Therefore,   
$[a,d]=a^{2^{n-1}}$ and
$(a^{i+2^{n-2}}\cdot d\m1)^2=b^2$.
 
If we replace $d$  by $d'=a^{i+2^{n-2}}d\m1$, we obtain 
$W=\gp{a,b,d'}$ and $(d')^2=b^2$, which is  case 2 
and we have that $W$ 
satisfies (\rm{iii}) of the  Theorem.

This proves that  the subgroup $W$ has a
generator system $u,v,w$ such that
$$
W=\gp{w,u,v\mid w^4=1, w^2=v^2, v^w=v\m1, u^{2^n}=1, 
 u^w=u^{1+2^{n-1}}, [u,v]=1},
$$
with $n>2$ 
and $H=\gp{u,v}$. 

Suppose that there exists $g\in G\setminus W$. Clearly, 
$G^{\prime}\subseteq W$ and $W$ is normal in $G$. Above we 
proved that the squares of all elements of $G$ outside $W$ 
belong to $H$ and they are central in $W$. Therefore, 
by the above argument we get  that $g^2=u^{2s}v^{2t}$  
for every $g\in G\setminus W$,
where $t\in\{0,1\}$, $s\in \Bbb N$. 
It is easy to see that 
$$
(g\m1u^s)^2=[g,u^s]g^{-2}u^{2s}=[g,u^s]v^{-2t}.
$$
If $(g\m1u^s)^2=1$ then
$g\m1u^s\in \Omega(W)\subseteq W$ and $g\in W$, which is impossible. 
Clearly, the order of elements $[g,u^s]$ and $v^{-2t}$ divide $2$ 
and $g\m1u^s$ is an element of order $4$. Then 
$M=\gp{g\m1u^s,v,w, u^{2^{n-2}}}$ is a subgroup of exponent $4$ 
and $\Omega(M)=\Omega(G)$. Clearly, $M/\Omega(M)$ 
is an  elementary abelian $2$-subgroup of order $16$. 
Therefore,  $M$ is a group with four generators. By O'Brien's Lemma,  
$M$ is isomorphic either to  $Q_8\times Q_8$  or to 
$H_{245}$ or to the central product $S_{2,2}$ and  $Q_8$ with common  
$a^2b^2$. It is impossible, because the centres 
of these groups 
have exponent $2$ but  in $M$ there
exists a central element $u^{2^{n-2}}$ of order $4$. It completes the
description of finite good  groups.

Now we suppose that $G$ is an infinite good group. We shall prove
that $G$ is the direct product of the quaternion group of order $8$
and the quasicyclic $2$-group.

It is easy to see that if $G$ has exponent 4 then $G$ is finite.
Indeed,  if $G$ is abelian then its finiteness follows from the
first Pr\"ufer's Theorem  and the condition 
$|\Omega(G)|\leq 4$. If
$G$  is nonabelian then take an ascending chain $G_1\subset
G_2\subset \cdots  $ of finite subgroups of $G$. It follows from
the description of finite good  groups given above that the
above chain is finite,  that is $G_n=G_{n+1}=\cdots$ for some $n\in
\Bbb N$. Since  $G$ is  locally finite,  $G=G_n$,  a contradiction.

Thus,  we may suppose that $\exp{(G)}>4$. We show now that
$\zeta(G)$ contains  a divisible subgroup. Let
$\mho_2(G)=\gp{g\in G\mid g^4=1}$ and 
$\Omega_2(G)=\gp{g^4\mid
g\in G}$.  Then $(gh)^4=g^4h^4$ for all $g, h\in G$ and the map
$g\to g^4$ is a group homomorphism of $G$ onto 
$\Omega_2(G)$ with kernel $\mho_2(G)$. Since $\Omega_2(G)$ has exponent 
$4$,  by the
above paragraph,  $\Omega_2(G)$ is finite. If   $\zeta(G)$ does not
contain a divisible group,  $\zeta(G)$ is finite. Hence, 
$\mho_2(G)\subseteq \zeta(G)$ is finite too which implies the
finiteness of $G$,  a contradiction.

We have that $\zeta(G)=R\times P$,  where $1\not=P$ is divisible and
$R$ does not contain a divisible subgroup. Observe now that
$R\not=1$ is cyclic,  $P$ is quasicyclic and for every noncentral
$g\in G$ there exists a noncentral element $g_1\in G$ such that 
$$
g_1\equiv g\pmod P  \text{  and   } \gp{g_1^2}=R.\tag3 
$$ 
Indeed,  we have that $g^2=cc'$,  where $c\in R,  c'\in P$.  Taking
$g_1=gd\m1$ with $d^2=c'$ we get $g_1^2=c$,  $g_1\equiv g \pmod P$.
As $g_1$ is noncentral, $c\not=1$. It follows that  
$R\not=1$ and since  $|\Omega(G)|=4$,  $R$ is cyclic and
$P$ is a quasicyclic. 
Let 
$R=\gp{z\mid z^{2^n}=1}$ with $n\geq 1$
and $P=\gp{c_1, c_2, \cdots, c_k, \cdots\mid c_1^2=1,  c_{k+1}^2=c_k}$
with $k=1, 2, \cdots$.
If  $c=z^l$ with  even $l$,  then $g_1z^{-\frac{l}{2}}$ is a
noncentral element of order $2$,  which is impossible. Hence,  $l$ is
odd and $\gp{g_1^2}=R$ as desired in   (3).

Next we observe that $R$ has order $2$ (that is $n=1$). Let $g$ and
$t$ be two noncommuting elements in $G$. We have that
$[g, t]=z^{2^{n-1}i}c_1^j$ ($i, j\in
\{0, 1\}$) and by   (3) we can suppose that $g^2=z^{i_1}$,  
$t^2=z^{i_2}$ with $i_1\equiv
i_2\equiv 1\pmod 2$.  Choose $m_1$ and $m_2$,  such that
$i_1m_1\equiv i_2m_2\equiv 1\pmod{2^n}$. If $n>1$,  then
$x=g^{m_1}t^{-m_2}z^{i2^{n-2}}c_2^j$ has order $2$ and
$[x, t]\not=1$,  a contradiction. 
Thus $R$ has order $2$ and $g^2=t^2=z$.
If $i=0$ then $gt\m1c_2^j$ is a noncentral 
element of order $2$, a contradiction.   Therefore, $[g, t]=zc_1^j$ 
and $Q=\gp{gc_2^j,  tc_2^j}$ is
isomorphic to the quaternion group of order $8$.

We show now that $G=Q\times P$ and it will complete  the proof of
necessity of the Theorem. Fix a noncentral element $x$ of $G$. 
There exists an element $c\in P$ such that $(xc)^2=(gc_2^j)^2$. Indeed,
if $j=0$ this  follows from   (3). In case $j=1$ by   (3) take $xc^{\prime}$ 
such that $(xc^{\prime})^2=z$. Then $(xc^{\prime}c_2)^2=zc_1$ as we need.

It is enough to show
that $x_1=xc\in W=Q\times \gp{c_3}$. Suppose that $x_1\not\in W$.
Then $W_1=\gp{x_1, W}$ has order $128$. Since  $W_1$ contains an element
of order $8$,  by the description of the finite good  groups given
above,  $W_1$ is one of the finite groups given in \rm{(ii)} and
\rm{(iii)} of the Theorem.
 
If $W_1=Q_1\times \gp{d\mid d^{16}=1}$,  where $Q_1$ is isomorphic
to the quaternion group of order $8$,  then $W=Q_1\times\gp{d^2}$
and as $x_1\not\in W$,  $x_1=vd$ for some $v\in W$. Hence, 
$x_1^8=d^8\not=1$,  contradicting the fact that $x_1^4=1$.

If $W_1$ is given by \rm{(iii)} of the Theorem,  then
$W=\gp{a, b}\times \gp{d^2}$ and $x_1=vd$ for some $v\in W$. Then we
have again that $x_1^8=d^8\not=1$ which is impossible.

Thus,  $x_1\in W$ and,  consequently,  $G\cong Q\times C_{2^\infty}$.

\subhead
3. Proof of the `if' part of the theorem
\endsubhead
We shall prove that $\Omega(V)=1+{\frak I}\big(\Omega(G)\big)$, 
where ${\frak I}\big(\Omega(G)\big)$ 
is the ideal generated by all elements of  form $g-1$ 
with $g\in \Omega(G)$.   
If $|\Omega(G)|=2$ then ${\frak I}^2\big(\Omega(G)\big)=0$ and 
$1+{\frak I}\big(\Omega(G)\big)\subseteq \Omega(V)$. 

Let $\Omega(G)=\gp{c\mid c^2=1}\times \gp{d\mid d^2=1}$. 
Then any $x\in {\frak I}\big(\Omega(G)\big)$
can be written as $x=\alpha(c+1)+\beta(d+1)+\gamma(c+1)(d+1)$, where
$\alpha,\beta,\gamma\in KG$. It is easy to see that 
$x^2=(\alpha\beta-\beta\alpha)(c+1)(d+1)=0$
and $1+{\frak I}\big(\Omega(G)\big)\subseteq \Omega(V)$.

\proclaim
{Lemma 4}
Let $KG$ be the   group algebra of an abelian $2$-group $G$ 
over a field  $K$ of characteristic $2$. 
If $x\in KG$, then $x^2=0$ if and only if 
$x\in {\frak I}\big(\Omega(G)\big)$. 
Moreover,
$$
\Omega(V)=1+{\frak I}\big(\Omega(G)\big).
$$ 
\endproclaim 

\rightline{\text{\qed}}

\proclaim
{Lemma 5} 
Let $H=\gp{c}$ be cyclic of order $2^n$ and let 
$v^2\in vKH(1+c^{2^{n-1}})$ for some $v\in KH$. Then $v^2=0$.
\endproclaim

\proof 
Let us  write $v$ as $\sum_{i=1}^{2^{n-1}-1}v_ic^i$ and $v_i\in 
K\gp{c^{2^{n-1}}}$. By induction on $k$ $(2\leq k\leq n)$ 
we shall  prove that 
$$
v(1+c^{2^{n-1}})=\sum_{i=0}^{2^{n-k}-1}\chi(v_i)c^i(1+c^{2^{n-1}})
(1+c^{2^{n-2}})\cdots (1+c^{2^{n-k}}). \tag4  
$$
Clearly, $v_i^2=\chi(v_i)^2$ and $v^2=\sum_{i=1}^{2^{n-1}-1}\chi(v_i)^2c^{2i}$.
Since $v^2\in vKH(1+c^{2^{n-1}})$ we obtain  that 
$\chi(v_i)^2=\chi(v_{i+2^{n-2}})^2$  
 which  implies that $\chi(v_i)=\chi(v_{i+2^{n-2}})$, where 
$i=0, 1, \ldots, 2^{n-2}-1$. Therefore, 
$$
v(1+c^{2^{n-1}})=\sum_{i=0}^{2^{n-2}-1}\chi(v_i)c^i(1+c^{2^{n-1}})
(1+c^{2^{n-2}})
$$
and the equality   (4) is true for $k=2$. Assume that 
$$
v(1+c^{2^{n-1}})=\sum_{i=0}^{2^{n-k+1}-1}\chi(v_i)c^i(1+c^{2^{n-1}})
(1+c^{2^{n-2}})\cdots (1+c^{2^{n-k+1}}).\tag5
$$
Then
$
\chi(v_i)=\chi(v_{i+{2^{n-k+1}}})=\chi(v_{i+{2^{n-k+2}}})=
\ldots =\chi(v_{i+{2^{n-2}}}), 
$
where $i=0,\ldots , 2^{n-1}-2^{n-2}-1$
and using this equality we get  
$$
v^2=\sum_{i=1}^{2^{n-1}-1}\chi(v_i)^2c^{2i}=
\sum_{i=1}^{2^{n-k+1}-1}\chi(v_i)^2c^{2i}(1+c^{2^{n-1}})
(1+c^{2^{n-2}})\cdots (1+c^{2^{n-k+2}}). \tag6
$$
Since $v^2\in vKH(1+c^{2^{n-1}})$, by (5) we have 
$$
v^2\in KH(1+c^{2^{n-1}})(1+c^{2^{n-2}})\cdots (1+c^{2^{n-k+1}})
$$
and by (6)  $\chi(v_i)= \chi(v_{i+2^{n-k}})$ for 
$i=1,2,\ldots, 2^{n-k}-1$. It implies by (5)  that   (4) is proved for 
all $k$.

Let $k=n$. Then $v(1+c^{2^{n-1}})=\chi(v_0)(1+c+c^2+\ldots+c^{2^n-1})$ 
and $v^2$ belongs to the ideal $K(1+c+c^2+\ldots+c^{2^n-1})$. Clearly,
$v^2=\alpha(1+c+c^2+\ldots+c^{2^n-1})$ for some $\alpha\in K$. 
Since $v^2\in K\gp{c^2}$ we conclude that $\alpha=0$ and $v^2=0$.

\rightline{\text{\qed}}

We remark that for any $G$  from the Theorem, 
$G/\zeta(G)$ is elementary abelian. For any 
$x=\sum_{g\in G}\alpha_gg\in \Omega(V)$  
we have   $x^2=\sum_{g\in G}\alpha_g^2g^2+w=1$, 
where $w$ belongs to the commutator subspace ${\Cal L}(KG)$. 
Obviously,  $g^2\in \zeta(G)$ and 
$\sum_{g\in G}\alpha_g^2g^2\in K\zeta(G)$.
We know that ${\Cal L}(KG)\cap K\zeta(G)=0$. Therefore, 
$$
\sum_{g\in G}\alpha_g^2g^2=1 \text{   and    }  w=0.\tag7
$$

We have the following cases:

Case 1. Let $G=\gp{a,b\mid a^{2^n}=b^{2^m}=1, a^b=a^{1+2^{n-1}}}$ with  
$n,m\geq 2$. Any $x\in \Omega(V)$ can be written as 
$x=x_0+x_1a+x_2b+x_3ab$, where $x_i\in K\zeta(G)$, and by   (7) we have
$$
x^2=x_0^2+x_1^2a^2+x_2^2b^2+x_3^2a^{2+2^{n-1}}b^2=1
$$
and $x_i^2\in K\zeta^2(G)$. Since
$\zeta(G)=\zeta^2(G)\cup a^2\zeta^2(G)
\cup b^2\zeta^2(G)\cup a^2b^2\zeta^2(G)$
we have
$x_0^2=1$, $x_1^2=x_2^2=x_3^2=0$
and by Lemma 4  we conclude that $\Omega(V)=1+{\frak I}\big(\Omega(G)\big)$.

Case 2. Now let 
$G=\gp{a,b\mid a^4=1, b^2=a^2, a^b=a\m1}\times \gp{c\mid c^{2^n}=1}$
with $n\geq 0$.
Any $x\in\Omega(V)$ can be written as    
$$
x=x_0+x_1a+x_2b+x_3ab,
$$ 
where $x_i\in K\zeta(G)$. By   (7) we have
$x^2=x_0^2+(x_1^2+x_2^2+x_3^2)a^2=1$ and
$$
x_0^2=1,\,\,\,\,\, 
x_1^2+x_2^2+x_3^2=0,\,\,\,\,\,
x_ix_j(1+a^2)=0, 
$$
where $i,j=1,2,3$ and $i\not=j$. 
Then $\chi(x_0)=1$ and $\chi(x_1)+\chi(x_2)+\chi(x_3)=0$. 

We shall prove that $x_1^2=x_2^2=x_3^2=0$.
Clearly,  $\zeta(G)=\gp{a^2}\times \gp{c}$ and  
$y^2\in K\gp{c^2}$, for any $y\in K\zeta(G)$.
Since $x_ix_j(1+a^2)=0$ ($i\not=j$),
we obtain that  $x_ix_j=z_{ij}(1+a^2)$, where $z_{ij}\in K\zeta(G)$. 
Clearly, by Lemma 4, 
$
x_1+x_2+x_3\in {\frak I}\big(\Omega(\zeta(G))\big)
$ and 
$$
x_1+x_2+x_3=z_1(c^{2^{n-1}}+1)+z_2(a^2+1),
$$
where $z_1,z_2\in K\gp{c}$.
Thus 
$$ 
x_1^2+x_1x_2+x_1x_3=z_1x_1(c^{2^{n-1}}+1)+z_2x_1(a^2+1) 
$$ 
and 
$x_1^2\equiv z_1x_1(c^{2^{n-1}}+1) \pmod {{\frak I}(\gp{a^2})}$, 
where ${\frak I}(\gp{a^2})=K\zeta(G)(1+a^2)$.   

Since $K\zeta(G)/{\frak I}(\gp{a^2})\cong K\gp{c}$, by Lemma 5  
$x_1^2\equiv 0 \pmod {{\frak I}(\gp{a^2})}$.
Clearly $x_i^2\in K\gp{c^2}$ and $K\gp{c^2}\cap  
{\frak I}(\gp{a^2})=0$. Thus, $x_1^2=0$ and similarly we conclude 
that  $x_i^2=0$, where $i=2, 3$.

Therefore, by Lemma 4 $\Omega(V)=1+{\frak I}\big(\Omega(G)\big)$.

Case 3.  Let $G=\gp{a,b, d \mid a^4=d^{2^n}=1, a^2=b^2, a^b=a\m1, 
d^a=d^{1+2^{n-1}}, b^d=b}$.  Any $x\in \Omega(V)$
 can be written as
$$
x=x_0+x_1a+x_2b+x_3ab +x_4d+x_5ad+x_6bd+x_7abd,
$$
where $x_i\in K\zeta(G)$. By   (7) we have
$$
x^2=x_0^2+(x_1^2+x_2^2+x_3^2)a^2
+x_4^2d^2+(x_5^2+x_7^2)a^2d^{2+2^{n-1}}+x_6^2d^2a^2=1\tag8
$$
and
$$
\cases
(x_1x_3+x_5x_7d^{2+2^{n-1}})(1+a^2)=0;\\
(x_1x_5+x_3x_7)(1+d^{2^{n-1}})=0;\\
(x_1x_7a^2+x_3x_5)(1+a^2d^{2^{n-1}})=0.
\endcases\tag9
$$
It is easy to see that 
$KG/{\frak I}(\gp{d^{2^{n-1}}})\cong  K(Q_8\times C_{2^{n-1}})$. By case 2, 
$x-1\in {\frak I}(\gp{a^2,d^{2^{n-2}}})$, and  
$$
(x_0+x_4d)+(x_1+x_5d)a+(x_2+x_6d)b+ (x_3+x_7d)ab\equiv 
1 \pmod {{\frak I}(\gp{a^2,d^{2^{n-2}}})}.
$$
We obtain that 
$$
\chi(x_0)+\chi(x_4)=1 \text{  and   }
\chi(x_i)+\chi(x_{i+4})=0\tag10 
$$
for all $i=1,2,3$.

First, suppose that $n=2$. Then $d^4=1$ and 
from   (8),   (9) and   (10) we have that either $\chi(x_i)=0$ 
($i=1,\cdots,7$) or $\chi(x_1)=\chi(x_3)=\chi(x_5)=\chi(x_7)$ 
and $\chi(x_2)=0$.

If  $\chi(x_1)=\chi(x_3)=\chi(x_5)=\chi(x_7)\not=0$ and $\chi(x_2)=0$, 
then $x_1,x_3,x_5,x_7$ are units  
and $x_i\m1=\chi(x_1)\m1x_i$ ($i=1,3,5,7$). Then    (9) implies  
$$
\cases
(1+\chi(x_1)\m1x_1x_3x_5x_7)(1+a^2)=0;\\
(1+\chi(x_1)\m1x_1x_3x_5x_7)(1+d^2)=0;\\
((1+a^2)+1+\chi(x_1)\m1x_1x_3x_5x_7)(1+a^2d^2)=0.
\endcases\tag11 
$$
Since $1+a^2d^2=(1+a^2)(1+d^2)+(1+a^2)+(1+d^2)$, by the third equality 
of   (11) $(1+a^2)(1+a^2d^2)=0$, which is impossible.

If now  $\chi(x_j)=0$ ($j=1,2,3,5,7$),  
then $1+x_0,x_j\in  {\frak I}\big(\Omega(G)\big)$ and 
by Lemma 4 
$x\in  1+{\frak I}\big(\Omega(G)\big)$.

Let $n>2$. 
Clearly,  $x\in 1+{\frak I}(\gp{a^2,d^{2^{n-2}}})$, 
$$
{\frak I}(\gp{a^2,d^{2^{n-2}}})=KG(a^2+1)+KG(1+d^{2^{n-1}})+KG(1+d^{2^{n-2}})
$$
and 
$$
x_i\equiv \sum_{j=0}^{2^{n-2}-1}\alpha_{ij}d^{2j}\;
(\mod {\frak I}(\gp{a^2,d^{2^{n-1}}})),
$$
where $\alpha_{ij}\in K$. Since $x\in 1+{\frak I}(\gp{a^2,d^{2^{n-2}}})$ we
have $\alpha_{ij}=\alpha_{i j+2^{n-3}}$, where  $i=1,\ldots, 7$, 
$j=0,\ldots, 2^{n-3}-1$. Thus $x_i$ can be  written as 
$$
x_i=u_i(1+a^2)+v_i(1+d^{2^{n-1}})+t_i(1+d^{2^{n-2}}),
$$
where $u_i,v_i\in K\zeta(G)$ and 
$$
t_i=\sum_{j=0}^{2^{n-3}-1}\alpha_{ij}d^{2j}.\tag12
$$
By   (8) 
$
x_0^2-1=x_1^2+x_2^2+x_3^2= x_4^2=(x_5^2+x_7^2)d^{2^{n-1}}+x_6^2=0
$
and by Lemma 4  
$x_0-1,x_4\in  {\frak I}(\gp{a^2,d^{2^{n-1}}})$. Thus 
$$ 
t_4^2(1+d^{2^{n-1}})=(t_1^2+t_2^2+t_3^2)(1+d^{2^{n-1}})=
(t_5^2+t_6^2+t_7^2)(1+d^{2^{n-1}})=0. \tag13
$$
Obviously,  $Supp(t_i^2)\subseteq \gp{d^4}$ and by   (13) we 
conclude that 
$t_4^2=t_1^2+t_2^2+t_3^2=t_5^2+t_6^2+t_7^2 =0$.
It follows that
$$
t_4=t_1+t_2+t_3=t_5+t_6+t_7=0.  \tag14
$$
We consider the projection of the equality $x^2=1$ onto cosets 
$a\zeta(G)$, $b\zeta(G)$, 
$ab\zeta(G)$, $d\zeta(G)$, 
$ad\zeta(G)$, $bd\zeta(G)$, $abd\zeta(G)$  
of $\zeta(G)$ in $G$.
Direct calculations show  
$$
\cases
(v_2t_3+v_3t_2+(u_4t_5+(u_6+v_6)t_7+(u_7+v_7)t_6)d^2)\times\\
\times(1+a^2)(1+d^{2^{n-1}})(1+d^{2^{n-2}})+
(t_2t_3+t_6t_7d^2)(1+a^2)(1+d^{2^{n-1}})=0;\\
(v_1t_3+t_1v_3+(v_5t_7+v_7t_5)d^2)\times\\
\times(1+a^2)(1+d^{2^{n-1}})(1+d^{2^{n-2}})+
(t_1t_3+t_5t_7d^2)(1+a^2)(1+d^{2^{n-1}})=0;\\
(v_1t_2+v_2t_1+(u_4t_7+(u_5+v_5)t_6+(u_6+v_6)t_5)d^2)\times\\
\times(1+a^2)(1+d^{2^{n-1}})(1+d^{2^{n-2}})+
(t_1t_2+t_5t_6d^2)(1+a^2)(1+d^{2^{n-1}})=0;\\
(u_1t_5+u_5t_1+u_3t_7+u_7t_3)(1+a^2)(1+d^{2^{n-1}})(1+d^{2^{n-2}})=0;\\
(u_4t_1+v_2t_7+v_7t_2+(u_3+v_3)t_6+(u_6+v_6)t_3)\times\\
\times(1+a^2)(1+d^{2^{n-1}})(1+d^{2^{n-2}})+
(t_2t_7+t_3t_6)(1+a^2)(1+d^{2^{n-1}})=0;\\
((u_1+v_1)t_7+(v_7+u_7)t_1+(u_3+v_3)t_5+(u_5+v_5)t_3)\times\\
\times(1+a^2)(1+d^{2^{n-1}})(1+d^{2^{n-2}})+
(t_1t_7+t_3t_5)(1+a^2)(1+d^{2^{n-1}})=0;\\
((v_1+u_1)t_6+(v_6+u_6)t_1+v_2t_5+v_5t_2+u_4t_3)\times\\
\times(1+a^2)(1+d^{2^{n-1}})(1+d^{2^{n-2}})+
(t_1t_6+t_2t_5)(1+a^2)(1+d^{2^{n-1}})=0.
\endcases\tag15
$$

Put $\omega_0=1$, $\omega_k=\omega_{k-1}(1+d^{2^{n-k-2}})$ $(k>0)$ and
$t_i^{(k)}=\sum_{j=0}^{2^{n-k-3}}\alpha_{ij}d^{2j}$, where  $k\geq 0$ and
the $\alpha_{ij}$  are the same as in   (12).

We shall establish by induction on $k$ that $t_i=t_i^{(k)}\omega_k$. 
Clearly, $t^{(0)}_i=t_i$. Using the definition of $\omega_k$ we have 
$$
(1+d^{2^{n-1}})(1+d^{2^{n-2}})\omega_{k-1}=
(1+d^{2^{n-1}})(1+d^{2^{n-2}})\cdots (1+d^{2^{n-k-1}}).
$$ 
By the induction hypothesis the first three equations of (15) 
imply 
$$
\cases
(t_2^{(k-1)}t_3^{(k-1)}  +t_6^{(k-1)}t_7^{(k-1)}d^2) &(1+d^{2^{n-1}})
\cdots (1+d^{2^{n-k-1}})=0;\\
(t_1^{(k-1)}t_3^{(k-1)}  +t_5^{(k-1)}t_7^{(k-1)}d^2) &(1+d^{2^{n-1}})
\cdots (1+d^{2^{n-k-1}})=0;\\
(t_1^{(k-1)}t_2^{(k-1)}  +t_5^{(k-1)}t_6^{(k-1)}d^2) &(1+d^{2^{n-1}})
\cdots (1+d^{2^{n-k-1}})=0.
\endcases\tag16  
$$
Because of   (14) we can write 
$$
(t_1^{(k-1)}+t_2^{(k-1)}+t_3^{(k-1)})\omega_{k-1}=
(t_5^{(k-1)}+t_6^{(k-1)}+t_7^{(k-1)})\omega_{k-1}=0,
$$
from which we obtain 
$
t_1^{(k-1)}+t_2^{(k-1)}+t_3^{(k-1)}=
t_5^{(k-1)}+t_6^{(k-1)}+t_7^{(k-1)}.
$
Using this equality we can put the expression of $t_3^{(k-1)}$ and 
$t_7^{(k-1)}$ into   (16), and obtain
$$
\cases
(t_2^{(k-1)}+t_6^{(k-1)}d)^2 &(1+d^{2^{n-1}})\cdots (1+d^{2^{n-k-1}})=0;\\
(t_1^{(k-1)}+t_5^{(k-1)}d)^2 &(1+d^{2^{n-1}})\cdots (1+d^{2^{n-k-1}})=0.\\
\endcases\tag17 
$$

The supports of the elements $(t_i^{(k-1)})^2$ and $d^2(t_j^{(k-1)})^2$ 
belong to different cosets of $\gp{d^4}$ in $\gp{d^2}$. Thus, by the equality  
$$
(t_i^{(k-1)})^2(1+d^{2^{n-1}})\cdots (1+d^{2^{n-k-1}})=
$$$$
\sum_{j=0}^{2^{n-k-3}-1}(\alpha_{ij}+\alpha_{ij+{2^{n-k-3}}})d^{4j}
(1+d^{2^{n-1}})\cdots (1+d^{2^{n-k-1}}), 
$$
it follows from   (17) that $\alpha_{ij}=\alpha_{ij+2^{n-3-k}}$ for all
$i=1,2,5,6$ and $j=0,\cdots, 2^{n-k-3}-1$. 
Hence, 
$$
t_i^{(k-1)}=\big(\sum_{j=0}^{2^{n-k-3}-1}\alpha_{ij}d^{2j}\big)
(1+d^{2^{n-k-2}})=
t_i^{(k)}(1+d^{2^{n-k-2}}), 
$$
and since by induction $t_i=t_i^{(k-1)}\omega_{k-1}$, we conclude that
$t_i=t_i^{(k)}\omega_k$,  as desired. 

By putting $k=n-3$ we get 
$t_i=\alpha_{i0}(1+d^2)\cdots (1+d^{2^{n-3}})$.  
Then taking $k=n-3$ in    (17) we obtain that $t_i=t_{i+4}=\alpha_i(1+d^2)\cdots
(1+d^{2^{n-3}})$, where $\alpha_i=\alpha_{i0}$ and $i=1,2,3$.

Then from   (15) we get the following system
$$
\cases
\alpha_1\chi(u_4+u_6+v_2+v_6)+\alpha_2\chi(u_6+u_7+v_2+v_3+v_6+v_7)+
\alpha_2^2+\alpha_1\alpha_2=0;\\
\alpha_1\chi(v_1+v_3+v_5+v_7)+\alpha_2\chi(v_1+v_5)+
\alpha_1^2+\alpha_1\alpha_2=0;\\
\alpha_1\chi(u_4+u_6+v_2+v_6)+\alpha_2\chi(u_4+u_5+v_1+v_5)+
\alpha_1\alpha_2=0;\\
\alpha_1\chi(u_1+u_3+u_5+u_7)+\alpha_2\chi(u_3+u_7)=0;\\
\alpha_1\chi(u_4+u_6+v_2+v_6)+\alpha_2\chi(u_3+u_6+v_2+v_3+v_6+v_7)=0;\\
\alpha_1\chi(u_1+u_3+u_5+u_7+v_1+v_3+v_5+v_7)+\alpha_2\chi(u_1+u_5+v_1+v_5)=0;\\
\alpha_1\chi(u_4+u_6+v_2+v_6)+\alpha_2\chi(u_1+u_4+v_1+v_5)=0.
\endcases
$$

It is easy to verify that the last  system yields
$\alpha_1=\alpha_2=0$. Therefore, $x\in 1+{\frak I}\big(\Omega(G)\big)$.

Case 4. Let 
$$
\eqalign{G=\gp{\,x,y,u\bigm|\ &x^4=y^4=1, x^2=[y,x],\cr
&y^2=u^2=[u,x], x^2y^2=[u,y]\,}.}
$$
It is easy to see that   $\zeta(G)=\gp{x^2, y^2}$.
Then any $\alpha \in \Omega(V)$ can be written as
$$
\alpha=\alpha_0+\alpha_1x+\alpha_2y+\alpha_3xy+\alpha_4u+
\alpha_5xu+\alpha_6yu+\alpha_7xyu,
$$
where $\alpha_i\in K\zeta(G)$. Clearly   $\alpha_i^2\in K$
 and by   (7) we have
$$
\alpha^2=\alpha_0^2+(\alpha_1^2+\alpha_5^2+\alpha_7^2)x^2+
(\alpha_2^2+\alpha_3^2+\alpha_4^2)y^2+\alpha_6^2x^2y^2=1\tag18 
$$
and
$$
\cases
\alpha_1\alpha_2(1+x^2)+\alpha_4\alpha_7(1+x^2)y^2=0;\\
\alpha_1\alpha_3(1+x^2)+\alpha_4\alpha_6(x^2+y^2)+
\alpha_5\alpha_7(1+y^2)=0;\\
\alpha_1\alpha_4(1+y^2)+\alpha_2\alpha_7(1+y^2)x^2=0;\\
\alpha_1\alpha_5(1+y^2)x^2+\alpha_2\alpha_6(x^2+y^2)+
\alpha_3\alpha_7(1+x^2)y^2=0;\\
\alpha_1\alpha_6(1+x^2y^2)+\alpha_2\alpha_5(1+y^2)x^2+
\alpha_3\alpha_4(1+x^2)=0;\\
\alpha_1\alpha_7(x^2+y^2)+\alpha_2\alpha_4(1+x^2y^2)=0;\\
\alpha_2\alpha_3(1+x^2)y^2+\alpha_4\alpha_5(1+y^2)+\alpha_6
\alpha_7(1+x^2y^2)=0.\\
\endcases\tag19 
$$
From   (18) we have 
$\alpha_1^2+\alpha_5^2+\alpha_7^2=
\alpha_2^2+\alpha_3^2+\alpha_4^2=\alpha_6^2=0$ and 
$\alpha_0^2=1$. Therefore,  
$\chi(\alpha_5)=\chi(\alpha_1)+\chi(\alpha_7)$
and $\chi(\alpha_3)=\chi(\alpha_2)+\chi(\alpha_4)$. 
If we multiply (19) first by  $1+x^2$ and then by $1+y^2$ 
and replace $\chi(\alpha_5)$ by $\chi(\alpha_1)+\chi(\alpha_7)$ 
and $\chi(\alpha_3)$ by $\chi(\alpha_2)+\chi(\alpha_4)$     
we obtain  $\chi(\alpha_j)=0$
for all $j=1,\cdots,7$.

Therefore,  $\Omega(V)=1+{\frak I}\big(\Omega(G)\big)$, which completes the
proof of the Theorem.

\rightline{\text{\qed}}

\enddocument

\address
\hskip-\parindent
Institute of Mathematics and Informatics
\newline
Lajos Kossuth University
\newline
H-410 Debrecen, P.O. Box 12
\newline
Hungary
\newline
vbovdi\@math.klte.hu
\endaddress
\address
\hskip-\parindent
Departamento de Matem\'atica
\newline
Universidade de S\~ao Paulo
\newline
Caixa Postal 66281,
\newline
CEP 05315--970, S\~ao Paulo
\newline
Brasil
\hskip-\parindent
\newline
dokucha\@ime.usp.br
\endaddress
\enddocument